\documentclass[a4paper]{amsart}

\usepackage{amssymb,latexsym}
\usepackage{amsmath}
\usepackage{amsfonts}
\usepackage{amsthm}
\usepackage{amssymb}

\theoremstyle{plain}
\newtheorem{theorem}{Theorem}

\newtheorem{proposition}[theorem]{Proposition}

\theoremstyle{definition}

\title{Doubly regular Diophantine quadruples}
\begin{document}

\date{}


\author[A. Dujella]{Andrej Dujella}
\address{
Department of Mathematics\\
Faculty of Science\\
University of Zagreb\\
Bijeni{\v c}ka cesta 30, 10000 Zagreb, Croatia
}
\email[A. Dujella]{duje@math.hr}

\author[V. Petri\v{c}evi\'c]{Vinko Petri\v{c}evi\'c}
\address{
Department of Mathematics\\
Faculty of Science\\
University of Zagreb\\
Bijeni{\v c}ka cesta 30, 10000 Zagreb, Croatia
}
\email[V. Petri\v{c}evi\'c]{vpetrice@math.hr}

\begin{abstract}
For a nonzero integer $n$, a set of $m$ distinct nonzero integers $\{a_1,a_2,\ldots,a_m\}$
such that  $a_ia_j+n$ is a perfect square for all $1\leq i<j\leq m$, is called a $D(n)$-$m$-tuple.
In this paper, by using properties of so-called regular Diophantine $m$-tuples
and certain family of elliptic curves,
we show that there are infinitely many essentially different sets
consisting of perfect squares which are simultaneously $D(n_1)$-quadruples
and $D(n_2)$-quadruples with distinct nonzero squares $n_1$ and $n_2$.
\end{abstract}


\maketitle

\section{Introduction}

For a nonzero integer $n$, a set of distinct nonzero integers $\{a_1,a_2,\ldots,a_m\}$ such that
$a_ia_j+n$ is a perfect square for all $1\leq i<j\leq m$, is called a  Diophantine $m$-tuple
with the property $D(n)$ or $D(n)$-$m$-tuple. Sometimes it is convenient to allow that $n=0$
in this definition.
The $D(1)$-$m$-tuples are called simply Diophantine $m$-tuples,
and sets of nonzero rationals with the same property are called rational Diophantine $m$-tuples.
The first rational Diophantine quadruple, the set $\left \{\frac{1}{16},\frac{33}{16},\frac{17}{4},\frac{105}{16}\right \}$,
was found by Diophantus of Alexandria.
By multiplying elements of this set by $16$ we obtain the $D(256)$-quadruple
$\{1,33,68,105\}$. The first Diophantine quadruple, the set $\{1,3,8,120\}$,
was found by Fermat. In 1969, Baker and Davenport \cite{BD69},
proved that Fermat's set cannot be extended to a Diophantine quintuple.
Recently, He, Togb\'e and Ziegler proved that there are no Diophantine quintuples~\cite{HTZ16+}.
Euler proved that there are infinitely many rational Diophantine quintuples.
The first example of a rational Diophantine sextuple, the
set $\{11/192, 35/192, 155/27, 512/27, 1235/48, 180873/16\}$, was found by Gibbs \cite{Gibbs1},
while Dujella, Kazalicki, Miki\'c and Szikszai \cite{DKMS} recently proved that there are infinitely
many rational Diophantine sextuples (see also \cite{Duje-Matija,DKP,DKP-split}). It is not known whether there
exists a rational Diophantine septuple.
Gibbs' example shows that there exists a $D(2985984)$-sextuple.
It is not known whether there exist a $D(n)$-septuple for some $n\neq 0$.
Moreover, it is not known whether there exist a $D(n)$-sextuple for any $n$ which is not a perfect square.
For an overview of results on
Diophantine $m$-tuples and its generalizations see \cite{Duje-Notices}.


In \cite{kk01}, A.~Kihel and O.~Kihel asked if there are Diophantine triples  $\{a, b, c\}$
which are $D(n)$-triples for several distinct $n$'s.
In \cite{ADKT}, several infinite families of Diophantine triples were presented
which are also  $D(n)$-sets for two additional $n$'s.
Furthermore, there are examples of Diophantine triples which are $D(n)$-sets for three additional $n$'s.
If we omit the condition that one of the $n$'s is equal to $1$,
then the size of a set $N$ for which there exists a triple $\{a, b, c\}$ of nonzero integers
which is a $D(n)$-set for all $n\in N$ can be arbitrarily large.

In \cite{DP-n1n2}, we proved that there are infinitely many nonequivalent sets of four distinct nonzero
integers $\{a, b, c, d\}$ with the property that there exist two distinct nonzero integers
$n_1$ and $n_2$ such that $\{a, b, c, d\}$ is a $D(n_1)$-quadruple and a $D(n_2)$-quadruple
(we called equivalent a quadruple $\{a,b,c,d\}$ with properties $D(n_1)$ and $D(n_2)$ and
a quadruple $\{au,bu,cu,du\}$ with properties $D(n_1u^2)$ and $D(n_2u^2)$ for a nonzero rational $u$).
We presented two constructions of infinite families of such quadruples.
The first of them contains pairs $\{a,b\}$ such that $a/b=-1/7$,
while in the second family we allowed that $n_1=0$.

In this paper, we will improve results of \cite{DP-n1n2} by considering so-called regular
Diophantine $m$-tuples. A (rational) $D(n)$-quadruple $\{a,b,c,d\}$ is called \emph{regular} if
\begin{equation} \label{regn}
n(d + c - a - b)^2 = 4(ab + n)(cd + n).
\end{equation}
Equation (\ref{regn}) is symmetric under permutations of $a, b, c, d$.
Since the right hand side of (\ref{regn}) is a square, it is clear that a regular
$D(n)$-quadruple may exist only if $n$ is a perfect square.
On the other hand, if $n=\ell^2$ is a perfect square, then e.g.
$\{\ell, 3\ell, 8\ell, 120\ell\}$ is a regular $D(\ell^2)$-quadruple.
A $D(\ell^2)$-quadruple $\{a,b,c,d\}$ is regular if and only if
the rational $D(1)$-quadruple $\{a/\ell, b/\ell, c/\ell, d/\ell\}$ is regular.

In this paper, we consider the question is it possible that a quadruple $\{a,b,c,d\}$
is simultaneously a regular $D(u^2)$-quadruple and a regular $D(v^2)$-quadruple for $u^2 \neq v^2$
(we called such sets \emph{doubly regular Diophantine quadruples}).
We will give an affirmative answer to this question.
Moreover, in our solution all elements $a,b,c,d$ will be perfect squares.
So, if we allow $n=0$ in the definition of $D(n)$-$m$-tuples,
we get quadruples which are simultaneously $D(n_1)$-quadruples, $D(n_2)$-quadruples
and $D(n_3)$-quadruples, with $n_1\neq n_2\neq n_3 \neq n_1$, thus improving the results
from \cite{DP-n1n2}.

Our main result is

\begin{theorem} \label{tm:1}
There are infinitely many nonequivalent sets of four distinct nonzero integers
$\{a,b,c,d\}$ which are regular $D(n_1)$ and $D(n_2)$-quadruples for
distinct nonzero squares $n_1$ and $n_2$. Moreover, we may take that all elements of
 these sets are perfect squares, so they are also $D(0)$-quadruples.
\end{theorem}

The construction of sets with the properties from Theorem \ref{tm:1}
use a parametrization of rational Diophantine triples
and properties of certain family of elliptic curves
(for other connections between Diophantine $m$-tuples and elliptic curves see e.g.
\cite{D-bordo,D-Peral}).

\section{Construction of doubly regular Diophantine quadruples}


As we mentioned in the introduction, in \cite{DP-n1n2} we constructed two
infinite families of such quadruples which are $D(n_1)$ and $D(n_2)$-quadruples with $n_1\neq n_2$.
We also listed some sporadic examples which do not fit in these two infinite families.
None of these examples is such that $n_1$ and $n_2$ are both nonzero squares.
However, in some of them one of the numbers $n_1$, $n_2$ is a square.
For example, $\{28, 6348, 18750, 88872\}$ is a $D(330625)$ and $D(38101225)$-quadruple
and $330625=575^2$. Moreover, $\{28, 6348, 18750, 88872\}$ is a regular
$D(330625)$-quadruple.

Assume now that $\{a_1,b_1,c_1,d_1\}$ is a regular $D(u^2)$-quadruple and regular $D(v^2)$-quadruple.
Then $\{a,b,c,d\}$, where $a=a_1/u$, $b=b_1/u$, $c=c_1/u$, $d=d_1/u$, is a regular
rational $D(1)$-quadruple, and $\{a/x,b/x,c/x,d/x\}$, where $x=v/u$, is also a regular
rational $D(1)$-quadruple.

We will use a parametrization of rational $D(1)$-triples which is a slight modification of
the parametrization due to L.~Lasi\'c \cite{luka} (see also \cite{DKP-split}).
Lasi\'c's parametrization is
\begin{align*}
a &=  \frac{2 t_1 (1 + t_1 t_2 (1 + t_2 t_3))} {(-1 + t_1 t_2 t_3) (1 + t_1 t_2 t_3)}, \\
b &=  \frac{2 t_2 (1 + t_2 t_3 (1 + t_3 t_1))} {(-1 + t_1 t_2 t_3) (1 + t_1 t_2 t_3)}, \\
c &=  \frac{2 t_3 (1 + t_3 t_1 (1 + t_1 t_2))} {(-1 + t_1 t_2 t_3) (1 + t_1 t_2 t_3)}.
\end{align*}
From the condition that $\{a,b,c,d\}$ is a regular $D(1)$-quadruple,
we compute $d$ and we obtain
$$ d=\frac{2(1+t_1t_2t_3)(t_1t_2+1+t_2)(t_1+1+t_3t_1)(1+t_3+t_2t_3)}{(-1+t_1t_2t_3)^3}. $$
By inserting these values of $a,b,c,d$ in the condition of regularity of
quadruple $\{a/x,b/x,c/x,d/x\}$, we obtain the following quartic equation in $x$:
\begin{equation} \label{quarticx}
 4x^4+(-a^2+2ab+2ad-b^2+2bc+2ac-c^2+2cd-d^2+2bd)x^2+4abcd = 0.
\end{equation}
By inserting the condition (\ref{regn}) with $n=1$ in the $x^2$-term in (\ref{quarticx}),
we obtain
\begin{equation*}
 4(x^2-1)(x^2-abcd) = 0.
\end{equation*}
Since we are interested in solutions with $u^2\neq v^2$, i.e. $x^2 \neq 1$,
we get that $x^2=abcd$. Thus, $abcd$ should be a perfect square, which leads to the condition that
$$
t_1t_2t_3(1+t_3+t_2t_3)(t_1+1+t_3t_1)(t_1t_2+1+t_2)(t_1t_2+t_1t_2^2t_3+1)(t_2t_3^2t_1+1+t_2t_3)(t_3t_1^2t_2+1+t_3t_1) $$
is a perfect square.

To solve the last condition, we introduce the following substitutions:
\begin{align*}
t_1 &= \frac{k}{t_2 t_3}, \\
t_2 &= m-\frac{1}{t_3}.
\end{align*}
Now the condition becomes
$$
k t_3 (1+m) (k+mt_3-1+kt_3) (k+t_3+mt_3-1) (km+1) (k+m) (k^2+mt_3-1+kt_3)= w^2, $$
which can be considered as a quartic in $t_3$:
\begin{equation} \label{quartic}
\begin{gathered}
k (m+1)^2 (km+1) (k+m)^3 t_3^4  \\
\mbox{}+ k (m+1) (k-1) (km+3m+2k+2) (km+1) (k+m)^2 t_3^3  \\
\mbox{}+ k (m+1) (k-1)^2 (km+1) (k+m) (k^2+2km+3k+3m+1) t_3^2  \\
\mbox{}+ k (m+1) (k+1) (k-1)^3 (k+m) (km+1) t_3 = w^2.
\end{gathered}
\end{equation}
The quartic (\ref{quartic}) has an obvious rational point $[t_3,w]=[0,0]$,
so it can be, in the standard way (see e.g. \cite[Section 1.2]{connell}),
transformed in an elliptic curve.
To ensure that this curve has positive rank,
we will force (\ref{quartic}) to have an additional rational point.
A good candidate for an additional point is $t_3=1/m$,
since it is a root of the discriminant of the left hand side of  (\ref{quartic})
with the respect to $k$.
By inserting $t_3=1/m$ in (\ref{quartic}), we get the condition that
$k(km+1)(k+m)$ is a perfect square (note that this condition is equivalent
to $ab$ being square). From $k(km+1)(k+m)=(km+z)^2$, we get
$m = \frac{k^2-z^2}{k(-1-k^2+2z)}$.
Here we take for the simplicity that $z=2$.

By transforming the quartic,
with the substitution
\begin{equation} \label{subst3}
 t_3 = k(k-1)(k+1)(k^2-3)(k^3-k^2-3k+4)/X, 
\end{equation} 
we obtain the following elliptic curve over $\mathbb{Q}(k)$:
\begin{align}
Y^2 &=
(X+ (k^3-k^2-3 k+4) (k^2-2)^2)
(X+ (k+1) (k^3-k^2-3 k+4) (k^2-2)^2) \nonumber \\
&\,\,\,\,\mbox{}\times (X+ (k+1) (k^3-k^2-3 k+4)^2) \label{YX}
\end{align}
with $2$-torsion points
\begin{align*}
T_1 &= [-(k+1)(k^3-k^2-3k+4)^2, 0], \\
T_2 &=[-(k+1)(k^3-k^2-3k+4)(k^2-2)^2, 0], \\
T_3 &=[-(k^3-k^2-3k+4)(-2+k^2)^2, 0],
\end{align*}
and an additional rational point
$$ P=[-(k-2) (k+2) (k+1) (k^3-k^2-3 k+4) (k-1), k^2 (k+1) (k^3-k^2-3 k+4)^2]. $$
The point $P$ does not give the desired solution because it corresponds
to $t_3=1/m$ which leads to $t_2=0$.
A point $[X,Y]$ would give us a solution if the corresponding quadruple $\{a,b,c,d\}$
satisfies that $ab+x^2, \ldots, cd+x^2$ are all perfect squares.
However, since $x^2=abcd$ and $ab+x^2=ab(cd+1)$, we see that the conditions are
equivalent to $ab,ac,ad,bc,bd,cd$ being perfect squares (i.e. to the condition that
$\{a,b,c,d\}$ is a $D(0)$-quadruple).
Since $ab=\frac{4(k^2-1)^2}{(k+1)^2(k-1)^2(k^2-3)^2}$ is a perfect square, 
and $ad = ac \cdot cd / c^2 = ac \cdot abcd /(c^2 \cdot ab)$, 
it suffices to satisfy the condition that $ac$ is a perfect square.
The condition is
$$
t_3 (4t_3-4t_3k^2-4k^3+3k+t_3k^4+k^5) = \Box, 
$$
which under substitution (\ref{subst3}) becomes 
$$ (k^3-k^2-3k+4) (X + (k^3-k^2-3k+4)(k^2-2)^2) = \Box. $$ 
Since this condition is satisfied for the $X$-coordinate of the point $P$, 
and $(X + (k^3-k^2-3k+4)(k^2-2)^2)$ is one of the factors of the right hand side of (\ref{YX}), 
by the $2$-descent argument (see \cite[Theorem 4.2]{Knapp}), 
it is satisfied for the values of $t_3$ which correspond to
$X$-coordinates of points of the form $P+2T$, hence it is satisfied for all 
odd multiples of the point $P$.

In particular, we may take the point
{\small
\begin{align*}
3P &=
\Big[\frac{1}{(k^6-6k^5-3k^4+28k^3-8k^2-32k+16)^2} \times (k-2) (k+2) (k-1) (k+1)  \\
& \,\,\,\,\mbox{}\times (3k^6-2k^5-13k^4+8k^3+16k^2-16) (5k^6-6k^5-27k^4+40k^3+32k^2-64k+16)\\
& \,\,\,\,\mbox{}\times (k^3-k^2-3 k+4), \\
& \frac{-1}{(k^6-6 k^5-3 k^4+28 k^3-8 k^2-32 k+16)^3} \times k^2 (k+1)  \\
& \,\,\,\,\mbox{}\times (4k^7-7k^6-22k^5+49k^4+20k^3-88k^2+32k+16) \\
& \,\,\,\,\mbox{}\times (4k^7-5k^6-26k^5+39k^4+48k^3-88k^2-16k+48) \\
& \,\,\,\,\mbox{}\times (k^6+2k^5-7k^4+8k^2-16k+16) (k^3-k^2-3k+4)^2
\Big]
\end{align*}}%
which corresponds to
{\small
$$ t_3=
\frac{k (k^2-3) (k^6-6k^5-3k^4+28k^3-8k^2-32k+16)^2} {(k-2) (k+2) (3k^6-2k^5-13k^4+8k^3+16k^2-16)
(5k^6-6k^5-27k^4+40k^3+32k^2-64k+16)}.
$$}%
By solving the quadratic equation in $x$, we obtain $x=x_1/x_2$, where
{\small
\begin{align*}
 x_1 &=
 (k^2-2) (k^6+2k^5-7k^4+8k^2-16k+16) (k^6-6k^5-3k^4+28k^3-8k^2-32k+16) \\
& \,\,\,\,\mbox{}\times (4k^7-5k^6-26k^5+39k^4+48k^3-88k^2-16k+48) \\
& \,\,\,\,\mbox{}\times (4k^7-7k^6-22k^5+49k^4+20k^3-88k^2+32k+16), \\
x_2 &= 2 (k+1) (k^2-3) (k^3-k^2-2k+4) (2k^4-k^3-7k^2+4k+4) (k-2)^2 (k+2)^2 (k-1)^3 \\
& \,\,\,\,\mbox{}\times (3k^6-2k^5-13k^4+8k^3+16k^2-16) (5k^6-6k^5-27k^4+40k^3+32k^2-64k+16) .
\end{align*}}%
By getting rid of denominators in $a,b,c,d,x$
we obtain the following proposition,
which clearly implies the statements of Theorem \ref{tm:1}.

\begin{proposition} \label{prop:1}
Let $k$ be an integer such that $k\neq 0, \pm 1, \pm 2$, and let \\
\begin{align*} a &=(k-1)^2 (k-2)^2 (k+2)^2(3k^6-2k^5-13k^4+8k^3+16k^2-16)^2 \\
& \,\,\,\, \mbox{} \times (5k^6-6k^5-27k^4+40k^3+32k^2-64k+16)^2, \\
b &=64 k^2 (k-1)^2 (k-2)^2 (k+2)^2 (k^3-k^2-3k+4)^2 (k^2-2)^2 \\
& \,\,\,\, \mbox{} \times (k^3-k^2-2k+4)^2 (2k^4-k^3-7k^2+4k+4)^2, \\
c &= k^2 (k-1)^2 (k^2-3)^2 (k^6-6k^5-3k^4+28k^3-8k^2-32k+16)^2 \\
& \,\,\,\, \mbox{} \times (4k^7-5k^6-26k^5+39k^4+48k^3-88k^2-16k+48)^2, \\
d &= (k+1)^2 (k^3-k^2-3k+4)^2 (k^6+2k^5-7k^4+8k^2-16k+16)^2 \\
& \,\,\,\, \mbox{} \times (4k^7-7k^6-22k^5+49k^4+20k^3-88k^2+32k+16)^2.
\end{align*}

Then $\{a,b,c,d\}$ is a $D(n_1)$, $D(n_2)$ and $D(n_3)$-quadruple, where
\begin{align*} n_1 &= 16k^2 (k+1)^2 (k-2)^4 (k+2)^4 (k-1)^6 (k^2-3)^2
\\ & \,\,\,\, \mbox{} \times
 (k^3-k^2-2k+4)^2 (k^3-k^2-3k+4)^2 (2k^4-k^3-7k^2+4k+4)^2
\\ & \,\,\,\, \mbox{} \times (3k^6-2k^5-13k^4+8k^3+16k^2-16)^2
\\ & \,\,\,\, \mbox{} \times (5k^6-6k^5-27k^4+40k^3+32k^2-64k+16)^2  , \\
n_2 &= 4k^2 (k^2-2)^2 (k^3-k^2-3k+4)^2 (k^6+2k^5-7k^4+8k^2-16k+16)^2
\\ & \,\,\,\, \mbox{} \times  (k^6-6k^5-3k^4+28k^3-8k^2-32k+16)^2
\\ & \,\,\,\, \mbox{} \times (4k^7-5k^6-26k^5+39k^4+48k^3-88k^2-16k+48)^2
\\ & \,\,\,\, \mbox{} \times (4k^7-7k^6-22k^5+49k^4+20k^3-88k^2+32k+16)^2  , \\
n_3 &= 0.
\end{align*}
\end{proposition}

For example, by taking $k=3$ in Proposition \ref{prop:1}, we obtain that
$$ \{1066758050, 7214407200, 8024417928, 44219811272\}$$
is a $D(90467582183447040000)$,
$D(30185892484109116209)$ and $D(0)$-quadruple.

\medskip

Other points $[X,Y]$ will not necessarily satisfy all required conditions.
However, for the point
{\small
\begin{align*}
P+T_1 &=
\Big[ -\frac{1}{(k^3-k^2-2k+4)^2} \times (k+1)(k^6+2k^5-7k^4+8k^2-16k+16) \\
& \,\,\,\,\mbox{}\times (k^3-k^2-3k+4)^2, \\
& -\frac{2}{(k^3-k^2-2k+4)^3} \times  k^2(k-2)(k+2)(k+1)(k^2-3) \\
& \,\,\,\,\mbox{}\times (2k^4-k^3-7k^2+4k+4)(k-1)^2(k^3-k^2-3k+4)^2/((k^3-k^2-2k+4)^3)
\Big]
\end{align*}}%
the corresponding $a,b,c,d,x$ satisfy that
$ab+x^2$, $cd+x^2$ are squares, while
$ac+x^2$, $ad+x^2$, $bc+x^2$ $bd+x^2$ are $(-k)\times {\rm squares}$.
By taking $k=-u^2$, we see that all conditions are satisfied,
and we obtain the following result.

\begin{proposition} \label{prop:2}
Let $u$ be an integer such that $u\neq 0, \pm 1$, and let \\
\begin{align*} a &= 2(u^6+2u^5+u^4-4u^2-4u-4)^2(u^6-2u^5+u^4-4u^2+4u-4)^2 \\
& \,\,\,\, \mbox{} \times (u^3-u^2+u-2)^2(u^3+u^2+u+2)^2, \\
b &= 2(2u^7-u^6+2u^5-u^4-6u^3+4u^2-8u+4)^2 \\
& \,\,\,\, \mbox{} \times (2u^7+u^6+2u^5+u^4-6u^3-4u^2-8u-4)^2 (u^4-2)^2, \\
c &= 2(u^2+1)^2(2u^8+u^6-7u^4-4u^2+4)^2 \\
& \,\,\,\, \mbox{} \times (u^6+u^4-2u^2-4)^2 u^2 (u^4-3)^2, \\
d &= 8(u-1)^2(u+1)^2u^2(u^4-3)^2(u^3-u^2+u-2)^2 \\
& \,\,\,\, \mbox{} \times (u^3+u^2+u+2)^2 (u^2+1)^4 (u^2+2)^2 (u^2-2)^2.
\end{align*}

Then $\{a,b,c,d\}$ is a $D(n_1)$, $D(n_2)$ and $D(n_3)$-quadruple, where
\begin{align*} n_1 &= (u-1)^2(u+1)^2(u^4-3)^2(u^2+1)^2(2u^7-u^6+2u^5-u^4-6u^3+4u^2-8u+4)^2 \\
& \,\,\,\, \mbox{} \times (2u^7+u^6+2u^5+u^4-6u^3-4u^2-8u-4)^2 (u^6+2u^5+u^4-4u^2-4u-4)^2 \\
& \,\,\,\, \mbox{} \times (u^6-2u^5+u^4-4u^2+4u-4)^2 (u^3-u^2+u-2)^2 (u^3+u^2+u+2)^2,  \\
n_2 &= 64(u^2+1)^4 (-2+u^4)^2 (2u^8+u^6-7u^4-4u^2+4)^2 (u^6+u^4-2u^2-4)^2 \\
& \,\,\,\, \mbox{} \times u^4(u^2+2)^2(u^2-2)^2 (u^4-3)^2 (u^3-u^2+u-2)^2 (u^3+u^2+u+2)^2, \\
n_3 &= 0.
\end{align*}
\end{proposition}

For example, by taking $u=2$ in Proposition \ref{prop:2}, we obtain that
$$ \{861184, 734247409, 15591268225, 8760960000\}$$
is a $D(30668429385921600)$,
$D(2816306908047360000)$ and $D(0)$-quadruple.


\medskip

Somewhat simpler examples can be found by a brute force search for
solutions $k,m,t_3$ of (\ref{quartic}) with small numerators and denominators.
Here are some examples obtained in that way:

\begin{center}
\begin{tabular}{r|l}
	$\{a,b,c,d\}$&$n_1,n_2,n_3$ \\		
	\hline
$\{1458, 66248, 5000, 14112\} $ & 16769025, 406425600, 0 \\
$\{451584, 25921, 12996, 950625\} $ & 30234254400, 4783105600, 0 \\
$\{985608, 11858, 57800, 352800\} $ & 49177497600, 4846248225, 0 \\
$\{105625, 50176, 72900, 1002001\} $ & 2981160000, 129859329600, 0 \\
$\{693889, 116964, 47089, 1982464\} $ & 144284503104, 52510639104, 0 \\
$\{74529, 2832489, 122500, 1115136\} $ & 134336910400, 214665422400, 0 \\
$\{438048, 3246152, 187272, 451250\} $ & 618173337600, 194388401025, 0 \\
$\{349448, 120050, 930248, 3645000\} $ & 493141017600, 288449555625, 0 \\
$\{31752, 45125000, 3426962, 18727200\} $ & 1409028350625, 65260546560000, 0 \\
$\{27766152, 1059968, 1820232, 61051250\} $ & 26694995558400, 122518001376225, 0
    \end{tabular}
	\end{center}

\bigskip

{\bf Acknowledgements.}
The authors want to thank to Matija Kazalicki and the referees for a careful reading of our paper
and for many valuable suggestions which improved the quality of the paper.
The authors were supported by the Croatian Science Foundation under the project no.~IP-2018-01-1313.
The authors acknowledge support from the QuantiXLie Center of Excellence, a project
co-financed by the Croatian Government and European Union through the
European Regional Development Fund - the Competitiveness and Cohesion
Operational Programme (Grant KK.01.1.1.01.0004).

\end{document}